\documentclass[conference]{ieeeconf}

\IEEEoverridecommandlockouts 

\usepackage{graphicx} 
\usepackage{epsfig} 
\usepackage{epstopdf}
\usepackage{booktabs} 
\usepackage{multirow} 
\usepackage{mathptmx} 
\usepackage{enumerate}
\usepackage{times} 
\usepackage{graphicx}
\usepackage{amsmath} 
\usepackage{amssymb} 
\usepackage{mathrsfs}
\usepackage{tikz}
\usepackage{latexsym}
\usepackage{pst-text}
\usepackage{dsfont}
\usepackage{pst-node}
\usepackage{pst-plot}
\usepackage{auto-pst-pdf}
\usepackage{amsthm}
\usepackage{caption}
\usepackage{float}
\usepackage{bm}
\usepackage{subfig}
\usepackage{verbatim}
\usepackage{cite}
\usepackage[T1]{fontenc}
\usetikzlibrary{shapes,arrows,chains,matrix,positioning,scopes,shadows,calc}
\newtheorem{theorem}{Theorem}

\newtheorem{lemma}{Lemma}

\newtheorem{corollary}{Corollary}
\newtheorem{definition}{Definition}
\newtheorem{assumption}{Assumption}

\usepackage[font=small,labelfont=bf]{caption}

\usepackage{calrsfs}
\DeclareMathAlphabet{\pazocal}{OMS}{zplm}{m}{n}

\usepackage{color}

\begin{document}

\title{Tracking of a Moving Target on a Hyperboloid and its Extension to Leader-follower Formation Controller}

\author{Ashish Shukla$^1$, Sanad Shaha$^{1\&2}$, Kashyap Rajeevsarathy$^3$, and Arijit Sen$^1$ %
\thanks{This project was funded by ANRF - Project No. ANRF/ARG/2025/007351/ENS.}
\thanks{$^{1}$ Department of Electrical Engineering and Computer Science, IISER Bhopal, India {\small ashish25, sanad22, ajsen@iiserb.ac.in}}%
\thanks{$^{2}$ TCS Research, India {\small sanad.shaha@tcs.com}}%
\thanks{$^{3}$ Department of Mathematics, IISER Bhopal, India {\small kashyap@iiserb.ac.in}}%
}

\pagestyle{empty}
\maketitle
\thispagestyle{empty}

\begin{abstract}
 This paper proposes a novel controller for a double-integrator agent to track a moving target with bounded velocity on a hyperboloid $\mathbb{H}^n$. The control law uses the intrinsic Riemannian structure of the hyperboloid, including the logarithm map, covariant derivative, and parallel transport, to make all its terms coordinate-free. Such formulation eases the implementation of the proposed tracking controller in the ambient Minkowski space without doing any additional chart computations. We use Barbalat's lemma for the hyperboloid to establish the convergence proof. Later, the proposed tracking controller is extended to a leader-follower formation control, where each agent reaches at a desired position around a moving leader in $\mathbb{H}^n$. We validate the results with numerical simulations considering two different target trajectories, geodesic and spiral, as examples.
\end{abstract}


\section{Introduction}\label{sec:prob_state:prelim}
Multi-agent systems (MAS) consist of multiple autonomous agents that interact and coordinate to achieve individual or collective objectives. Its applications span surveillance, search and rescue, environment monitoring, and autonomous vehicle coordination, and multi-robot containment \cite{savkin2021navigation, hua2025multi}. Target tracking and formation control constitute the core challenges in these applications.

The target tracking problems \cite{savkin2021navigation} refer to the pursuit of a moving object by continuously estimating its states, such as position, velocity, and acceleration, which are fed back to agents for generating the control inputs to intercept or follow the target. An alternative formulation is rendezvous with a moving target \cite{manathara2012rendezvous}, where, unlike the classical rendezvous at a static or mutually agreed point, agents continuously adapt their trajectories in response to the target's motion. When each agent is further assigned a prescribed offset from the target, the tracking problem transforms into a formation control problem \cite{besseghieur2020trajectory}, where agents must converge to their designated positions and preserve the geometric structure of the formation as the target moves.

While a major portion of the existing literature has addressed these problems in Euclidean space \cite{wang2022robust, vu2023distance}, there has been a growing interest in extending such formulations to Riemannian manifolds \cite{tron2012riemannian}. Recent works on target tracking and formation control in the Riemannian setting \cite{ahn2025quantitative, choi2023multi, li2024leader} either assume the flatness or compactness of the subspace of the manifold containing the target's trajectory. As per Bonnet-Myers theorem \cite{manfredo1992carmo}, every complete Riemannian manifold with uniformly positive sectional curvature is compact. Since flat manifolds have also been extensively studied \cite{vu2023distance, li2022leader}, it is important to investigate whether these results can be extended to negatively curved manifolds.

The hyperboloid $\mathbb{H}^n$ is a Riemannian manifold with a complete metric of constant sectional curvature $-1$. Because of its exponential volume growth, it has been extensively used in learning-based approaches \cite{ganea2018hyperbolic}. It also arises in modern physics, including celestial mechanics \cite{curtis2019orbital} and relativity \cite{ahn2023relativistic}. For certain aerospace missions involving cycler spacecraft, rendezvous and docking must be performed while the vehicles are traveling on escape or flyby trajectories relative to a planetary body \cite{landau2007guidance}. Such trajectories are hyperbolic Keplerian orbits and, in suitable coordinates, can be represented by points on a hyperboloid \cite{li2025rendezvous}. This connection highlights the potential relevance of cooperative control of MAS on the hyperboloid, motivating the study of cooperative tasks such as rendezvous, formation control, and target tracking in this setting.

Recent works on hyperboloids include rendezvous \cite{li2025rendezvous}, flocking \cite{ahn2021emergent, ahn2023relativistic}, and selective leader-follower convergence \cite{ahn2026selective} problems. To the best of our knowledge, the target-tracking problem for a target with an unbounded trajectory on hyperboloids has not yet been addressed. A recent work \cite{ahn2026selective} considers leader-follower coordination for possibly unbounded leader trajectories, where the control laws of leaders explicitly need feedback of followers' trajectories to guarantee convergence. Our controller eliminates such assumption. Because of the possibly unbounded trajectory of the moving target, LaSalle's invariance principle cannot be applied directly to these problems. To overcome this issue, we use the hyperboloid version of Barbalat's lemma \cite{li2025rendezvous}, which provides the necessary convergence guarantees in this non-compact setting. {\textit{Ahn} \cite{ahn2026selective} also extends Barbalat's lemma to address the selective leader-following problem using attractive and repulsive fields. This results in lack of explicit damping and consequently exhibits transient velocity oscillations. Explicit damping in our tracking controller eliminates transient oscillations as observed in \cite{ahn2026selective}.} 

The key features of this work are summarized as follows:
\begin{itemize}
\item To our knowledge, this is the first target-tracking controller for a double-integrator agent on $\mathbb{H}^n$, where the target may have a possibly unbounded trajectory with bounded velocity. The controller uses the logarithm map, parallel transport, and covariant derivative.

\item The controller is fully intrinsic and coordinate-free. Thus, it can be directly implemented in the ambient Minkowski space without any additional chart computations.

\item Convergence is guaranteed for any initial configuration of target and agent under an assumption that the target's velocity is bounded. This significantly relaxes the assumption in \cite{ahn2023asymptotic}, where the target trajectory is bounded.

\item We propose a virtual target-tracking protocol to extend this idea to leader-follower formation control, where the relative positioning of the virtual targets is preserved under parallel transport.
\end{itemize}
The rest of the paper is structured as follows: Section \ref{sec:prelim} presents the requisite mathematical notions on hyperbolic geometry used in this paper, which is followed by Section \ref{sec:problem formulation} where the problem formulation is presented. Section \ref{sec:main_results} proposes a nonlinear controller and proves its convergence. Section \ref{sec:leader-follower_formation} extends this controller to a leader-follower formation protocol. The performance of both proposed controllers is validated through simulations in Section~\ref{sec:sim}. Finally, Section~\ref{sec:conc} offers concluding remarks.
\section{Preliminaries on Hyperbolic Space}\label{sec:prelim}
This section presents some standard notation and results concerning the $n$-dimensional hyperbolic space used throughout the paper \cite{lee2018introduction}.
The $n$-dimensional hyperbolic space, denoted by $\mathbb{H}^n$, is the unique, simply connected, complete Riemannian manifold of constant sectional curvature $-1$. A standard realization of $\mathbb{H}^n$ is the hyperboloid model in the Minkowski space $\mathbb{R}^{n,1}$. The Minkowski space is an $(n+1)$-dimensional real vector space equipped with a symmetric bilinear form $\langle \cdot , \cdot \rangle _M$ defined as follows:
\begin{definition}\cite{lee2018introduction}
    Given two points $x=(x_0,x_1,\ldots,x_n)$ and $y=(y_0, y_1, \dots ,y_n)$ in $\mathbb{R}^{n,1}$, the Minkowski inner product $\langle x, y\rangle _M$ is defined as
    \begin{center}
        $\langle x, y\rangle _M = -x_0y_0+ \sum_{j=1}^n x_jy_j \in \mathbb{R}.$
    \end{center}    
\end{definition}
We identify $\mathbb{H}^n$ as the upper sheet of the two-sheeted hyperboloid in $\mathbb{R}^{n,1}$ given by
\begin{center}
    $\mathbb{H}^n= \{ x \in \mathbb{R}^{n,1} \vert \langle x , x \rangle_M = -1,~ x_0 >0\}.$
\end{center}
A vector $v \in \mathbb{R}^{n,1}$ is said to lie in the tangent space of $x \in \mathbb{H}^n$, denoted by $T_x\mathbb{H}^n$, if and only if $\langle x , v \rangle_M = 0$. The space $\mathbb{H}^n$ equipped with the metric $g_x(\cdot, \cdot):T_x\mathbb{H}^n \times T_x\mathbb{H}^n \rightarrow \mathbb{R}$ defined as $g_x(u,v) =\langle u, v\rangle _M$, becomes a Riemannian manifold. We denote the norm induced by $g$ as $\| \cdot \|_M$. {To describe motion and control intrinsically on $\mathbb{H}^n$, we require maps to relate points on the manifold to tangent vectors and vice versa. The following definitions introduce the same.}

\begin{definition}\cite{lee2018introduction}\label{def:log:exp}
    Let {$\gamma: [0,1] \to \mathbb{ H}^n$} be the geodesic from $x$ to $y$ in $\mathbb{H}^n$. The logarithm map evaluated at $x$, denoted by $\log_xy: \mathbb{H}^n \rightarrow T_x\mathbb{H}^n$, is defined as: $\log_xy := \dot{\gamma}(0).$ The exponential map is defined as the inverse map of the logarithm map. Thus, $\exp_x :T_x\mathbb{H}^n \rightarrow \mathbb{H}^n , ~\exp _x := \log_x ^{-1}$.
\end{definition}

\begin{definition}\cite{lee2018introduction}
    Given a smooth map $f: \mathbb{H}^n \rightarrow \mathbb{R}$, we define the gradient of $f$, denoted as $\operatorname{grad}f$, to be the unique vector field on $\mathbb{H}^n$ such that
    \begin{center}
        $\langle \operatorname{grad}f(x),v \rangle_M= \frac{d}{dt}f(\gamma(t)) \Big|_{t=0},$
    \end{center}
    $\forall~ x \in \mathbb{H}^n$ and $\forall~ v\in T_x\mathbb{H}^n$, where $\gamma:[0,1] \rightarrow\mathbb{H}^n$ is the geodesic such that $\gamma(0)=x$ and $\dot{\gamma}(0)=v$.
\end{definition}
{The following identity establishes a useful relationship between the Riemannian gradient of the squared distance function and the logarithm map \cite{fiorimanifold}:}
{\small \begin{equation}\label{eq:dist_grad}
    \begin{matrix}
        \frac{1}{2}\operatorname{grad}_x d(x,y)^2 = -\log_x y.
    \end{matrix}
\end{equation}}
Given two points $x$ and $y$ in $\mathbb{H}^n$, we denote the parallel transport from $x$ to $y$ as $P_x^y : T_x\mathbb{H}^n \rightarrow T_y\mathbb{H}^n$. This operator compares tangent vectors at two different points. {The operator $P_x^y$ is isometric, that is $\langle P_x^y u, P_x^y v\rangle_M = \langle u, v \rangle_M$ for all $u, v \in T_x \mathbb H^n$ \cite[Theorem 6]{fiorimanifold}. For proving the convergence results in this paper, we use the following extension of Barabalat's lemma on the hyperboloid.
\begin{lemma}\cite[Lemma 3.1]{li2025rendezvous}\label{lemma:barbalat:hyperboloid}
    Let $\{x_i(t), v_i(t)\}_{i=1}^N$ be a smooth curve on $\mathbb{H}^n$ with velocity $\dot{x}(t)=v(t)$ satisfying $\lim _{t \to \infty} v(t)=0$ and $\sup_{0 \leq t <\infty} \| \frac{D}{dt}\left(\frac{Dv}{dt} \right)\|_M < \infty$, then $\lim_{t \to \infty} \frac{Dv}{dt}=0$.
\end{lemma}}

\section{Target-tracking Problem Formulation}\label{sec:problem formulation}

Let $\pazocal{T}$ be a target moving in $\mathbb{H}^n$, with time-varying position $x_{\pazocal{T}}(t)\in \mathbb{H}^n$ and velocity $v_{\pazocal{T}}(t) \in T_{x_{\pazocal{T}}(t)}\mathbb{H}^n$ at time $t$. We introduce a double-integrator target-tracking protocol for an agent with position $x(t) \in \mathbb{H}^n$ and velocity $v(t) \in T_x\mathbb{H}^n$, with agent dynamics
{\small \begin{equation}\label{eq:double_int_agent}
\dot{x}(t) = v(t), \qquad     \frac{Dv}{dt}(t) = u(t),
\end{equation}}
where $u(t) \in {T_{x(t)}\mathbb{H}^n}$ is the control input and $D/dt$ is the covariant derivative. We define the position and velocity errors to track the target as follows:
{\small \begin{equation}\label{eq:error}
\tilde{x}(t) = \log_{x(t)} x_{\pazocal{T}}(t), \qquad \tilde{v}(t) = v(t) - P_{x_{\pazocal{T}}(t)}^{x(t)} v_{\pazocal{T}}(t).
\end{equation}}

\begin{assumption}\label{ass:target:speed:bound}
The target velocity $v_{\pazocal{T}}(t)$ is uniformly bounded. Thus, $\operatorname{sup}_{t\geq 0}\|v_{\pazocal{T}}(t)\|_M \leq M_v < \infty $. However, $x_{\pazocal{T}}(t)$ may have any value in $\mathbb{H}^n$.
\end{assumption}
Under Assumption \ref{ass:target:speed:bound}, we propose a control law such that an agent successfully tracks a target moving in $\mathbb{H}^n$:
{\small \begin{equation}\label{eq:tracking:defn}
    \lim _{t \rightarrow \infty} \|\tilde{x}(t)\|_M= \lim_{t \rightarrow \infty} \|\tilde{v}(t)\|_M=0.
\end{equation}}

\section{Target-tracking Algorithm in $\mathbb{H}^n$}\label{sec:main_results}
To achieve target tracking in $\mathbb{H}^n$, we propose the following nonlinear tracking control law:
{\small \begin{equation}\label{eq:controller}
    u(t) = \frac{D}{dt}\left(P_{x_{\pazocal{T}}(t)}^{x(t)} v_{\pazocal{T}}(t)\right) + k_p \log_{x(t)} x_{\pazocal{T}}(t) - k_d \tilde{v}(t),
\end{equation}}
where $k_p > 0$ and $k_d > 0$ are the control gains. {The controller \eqref{eq:controller} comprises three terms: the first term compensates for target acceleration, the second term is a proportional feedback that drives the agent toward the target, and the third is a damping term that regulates the velocity tracking error, suppressing oscillatory behavior.}
{
\begin{theorem} \label{thm:tracking}
Consider {an} agent following \eqref{eq:double_int_agent} and a moving target under Assumption \ref{ass:target:speed:bound}. Under the nonlinear controller \eqref{eq:controller}, the agent tracks the target on hyperboloid $\mathbb{H}^n$. Thus, 
{\small \begin{equation*}
    \lim _{t \rightarrow \infty} \|\tilde{x}(t)\|_M= \lim_{t \rightarrow \infty} \|\tilde{v}(t)\|_M=0.
\end{equation*}}
    
\end{theorem}
}

\begin{proof}
    To study the convergence under \eqref{eq:controller}, we introduce the following Lyapunov function {$V : T\mathbb{H}^n \times T\mathbb{H}^n \rightarrow \mathbb{R} $} for the agent-target system:
{\small \begin{equation}\label{eq:Lyapunov:target}
    V(t)=\frac{1}{2}k_pd^2(x(t), x_{\pazocal{T}}(t)) + \frac{1}{2} \|\tilde{v}(t) \|^2_M
\end{equation}}
 where $V(t)$ is the sum of a position error term and velocity error term. For convenience, we compute the derivative of each term separately to obtain $\dot{V}(t)$.

\noindent {\bf Derivative of position error term:} Using the identity \eqref{eq:dist_grad},
{\small    
 \begin{equation}\label{eq:position_error}
        \frac{d}{dt}\left(\frac{1}{2} d^2(x, x_{\pazocal{T}})\right) = \langle -\log_{x} x_{\pazocal{T}}, v \rangle_M + \langle -\log_{x_{\pazocal{T}}} x, v_{\pazocal{T}} \rangle_M.
    \end{equation}}
    For the second term of \eqref{eq:position_error}, using the invariance of parallel transport for the inner product on the tangent space, we get
{\small \begin{equation}\label{eq:position_error:int}
\begin{split}
            \langle -\log_{x_{\pazocal{T}}} x, v_{\pazocal{T}} \rangle_M= & \langle P_{x_{\pazocal{T}}(t)}^{x(t)}( -\log_{x_{\pazocal{T}}} x),~ P_{x_{\pazocal{T}}(t)}^{x(t)}(v_{\pazocal{T}}) \rangle_M\\
            =&\langle \log_{x} x_{\pazocal{T}},P_{x_{\pazocal{T}}(t)}^{x(t)}(v_{\pazocal{T}})\rangle_M.
        \end{split}
    \end{equation}}
Substituting \eqref{eq:position_error:int} into \eqref{eq:position_error}, we get
{\small  \begin{align}\label{eq:position_error_Final}
\frac{d}{dt}\left(\frac{1}{2} d^2(x, x_{\pazocal{T}})\right) = & \langle -\log_{x} x_{\pazocal{T}}, v \rangle_M +\langle \log_{x} x_{\pazocal{T}},P_{x_{\pazocal{T}}(t)}^{x(t)}(v_{\pazocal{T}})\rangle_M \nonumber\\
= & -\langle \log_{x} x_{\pazocal{T}},v-P_{x_{\pazocal{T}}(t)}^{x(t)}(v_{\pazocal{T}})\rangle_M \nonumber\\
\Rightarrow \frac{d}{dt}\left(\frac{1}{2} d^2(x, x_{\pazocal{T}})\right) = & -\langle \log_{x} x_{\pazocal{T}}, \tilde{v} \rangle_M
\end{align}}
\noindent \textbf{Derivative of the Velocity Error Term:} Since the Minkowski metric is compatible with the Levi-Civita connection, we take the covariant derivative directly inside the inner product and obtain the following:
{\small \begin{equation}\label{eq:velocity_error}
    \frac{d}{dt}\left(\frac{1}{2} \|\tilde{v}\|_M^2\right) = \left\langle \tilde{v}, \frac{D\tilde{v}}{dt} \right\rangle_M.
\end{equation}\label{eq:derivative_velocity_error}}
Expanding covariant derivative of velocity error, we get
{\small \begin{equation}\label{eq:derivative_velocity_error:cov}
\frac{D\tilde{v}}{dt} = \frac{D v}{dt} - \frac{D}{dt}\left(P_{x_{\pazocal{T}}}^{ x} v_{\pazocal{T}}\right) = u - \frac{D}{dt}\left(P_{x_{\pazocal{T}}}^{x} v_{\pazocal{T}}\right).
\end{equation}}
Substituting \eqref{eq:controller}, the derivative \eqref{eq:derivative_velocity_error:cov} becomes
{\small \begin{equation}\label{eq:derv:vtilde}
\begin{split}
\frac{D\tilde{v}}{dt} = &\left(\frac{D}{dt}\left(P_{x_{\pazocal{T}}}^{ x} v_{\pazocal{T}}\right) + k_p \log_{x} x_{\pazocal{T}} - k_d \tilde{v}\right) - \frac{D}{dt}\left(P_{x_{\pazocal{T}}}^{ x} v_{\pazocal{T}}\right) \\
= & k_p \log_{x} x_{\pazocal{T}} - k_d \tilde{v}.
\end{split}
\end{equation}}
Thus, the velocity error derivative evaluates to
{\small \begin{equation}\label{eq:velocity_error_final}
\begin{split}
\frac{d}{dt}\left(\frac{1}{2} \|\tilde{v}\|_M^2\right) = & \langle \tilde{v}, k_p \log_{x} x_{\pazocal{T}} - k_d \tilde{v}\rangle_M \\
= &k_p \langle \tilde{v}, \log_{x} x_{\pazocal{T}} \rangle_M - k_d \|\tilde{v}\|_M^2
\end{split}
\end{equation}}
Substituting \eqref{eq:position_error_Final} and \eqref{eq:velocity_error_final} into \eqref{eq:Lyapunov:target} gives
{\small \begin{equation*}
\begin{split}
\dot{V}(t)=&-k_p\langle \log_{x} x_{\pazocal{T}}, ~\tilde{v} \rangle+ k_p \langle \tilde{v}, ~\log_{x} x_{\pazocal{T}} \rangle_M - k_d \|\tilde{v}\|_M^2\\
= & - k_d\|\tilde{v}\|_M^2 \leq 0.
\end{split}
\end{equation*}}
This shows that $V(t)$ is non-increasing and therefore $V(t)\leq V(0),~\forall~ t\geq 0$. Differentiating $V$ twice, we get
{\small \begin{equation}\label{eq:barbalat:1}
\ddot{V}(t) = \frac{d}{dt}\left(-k_d \langle \tilde{v}, ~\tilde{v} \rangle_M\right) = -2k_d \left\langle \tilde{v},~ \frac{D\tilde{v}}{dt}\right\rangle_M
\end{equation}}
Substituting \eqref{eq:derv:vtilde} into \eqref{eq:barbalat:1} leads to
{\small \begin{equation}\label{eq:double_dot:before}
  \ddot{V}(t) = -2k_d k_p \langle \tilde{v}(t),~ \log_x x_{\pazocal{T}}(t) \rangle_M + 2k_d^2 \|\tilde{v}(t)\|_M^2.
\end{equation}}
Since $V$ is non-increasing, as $\dot{V}$ is negative semi-definite
\begin{equation*}
V(t)=\frac{1}{2}k_pd^2(x(t), x_{\pazocal{T}}(t)) + \frac{1}{2} \langle \tilde{v}(t), \tilde{v}(t) \rangle_M \leq V(0)
\end{equation*}
 Thus, $d(x(t),x_{\pazocal{T}}(t)) \leq D_0:=\sqrt{2V(0)/k_p}$ and $\| \tilde{v}(t) \|_M \leq \tilde{V} :=\sqrt{2V(0)}$. Now, applying the Cauchy-Schwarz inequality to the first term of \eqref{eq:double_dot:before} and using $\|\log_x x_{\pazocal{T}}\|_M=d(x,x_{\pazocal{T}})$ together, we get a uniform bound on $\ddot{V}(t)$ as
{\small \begin{equation}
\begin{split}
|\ddot{V}(t)| \leq & 2k_dk_p\|\tilde{v}(t)\|_Md(x(t),x_{\pazocal{T}}(t))+2k_d^2\|\tilde{v}(t)\|_M^2\\
\leq & 4k_d \sqrt{k_p}V(0)+4k^2_dV(0) < \infty.
\end{split}
\end{equation}}
Using Barbalat's lemma \cite{vidyasagar2002nonlinear}, we conclude that $\dot{V} \rightarrow 0$ as $t \rightarrow \infty$. Hence, $\lim_{t \rightarrow \infty}\| \tilde{v}(t)\|_M=0$. 

Now, we have to prove the convergence of $\lim_{t \rightarrow \infty} d(x(t),x_{\pazocal{T}}(t))=0$. From \eqref{eq:derv:vtilde}, we have
{\small \begin{equation*}
\begin{split}
\frac{D^2 \tilde{v}}{dt^2} = & \frac{D}{dt} (k_p \log_{x} x_{\pazocal{T}} - k_d \tilde{v})= k_p \frac{D}{dt}\left(\log_{x} x_{\pazocal{T}}\right) - k_d \frac{D\tilde{v}}{dt}\\
= & k_p \frac{D}{dt}\left(\log_{x} x_{\pazocal{T}}\right) - k_d k_p \log_{x} x_{\pazocal{T}} + k_d^2 \tilde{v}\\
\Rightarrow \left\| \frac{D^2\tilde{v}}{dt^2} \right\|_M \le & k_p \left\| \frac{D}{dt}\left(\log_{x} x_{\pazocal{T}}\right) \right\|_M + k_d k_p \|\log_{x} x_{\pazocal{T}}\|_M + k_d^2 \|\tilde{v}\|_M
\end{split}
\end{equation*}}
 From the previous part of the proof, we have $d(x,x_{\pazocal{T}})$ and $ \|\tilde{v}\|_M$ are bounded. Thus, $\left\| D^2\tilde{v}/dt^2\right\|_M $ becomes bounded if and only if $\left\| D\left(\log_{x} x_{\pazocal{T}}\right)/dt \right\|_M$ is bounded. The next part proves the same under Assumption \ref{ass:target:speed:bound}.

For brevity, denoting $d(x(t),x_{\pazocal{T}}(t))$ by $d$, we get from \eqref{eq:position_error_Final}
{\small \begin{equation} \label{eq:ddot_bound}
\begin{split}
 d \dot{d} =&-\langle \log_{x} x_{\pazocal{T}}, \tilde{v} \rangle \Rightarrow \dot{d} = -\langle \log_{x} x_{\pazocal{T}}, \tilde{v} \rangle/d\\
 \implies |\dot{d}| \leq & \|\tilde{v}\|_M \leq \tilde{V} =\sqrt{2V(0)}
 \end{split}
\end{equation}}
 Note that we used the Cauchy-Schwarz inequality and the fact $\|\log_x x_{\pazocal{T}}\|_M = d$ to get bounds of $|\dot{d}|$ in \eqref{eq:ddot_bound}. Furthermore, from the definition \eqref{eq:error}, $v = \tilde{v} + P^{x}_{x_{\pazocal{T}}} v_{\pazocal{T}}$, we get the following by applying the triangle inequality:
{\small \begin{equation}\label{eq:bound:v}
    \|v(t)\|_M \leq \|\tilde{v}(t)\|_M + \|v_{\pazocal{T}}(t)\|_M \leq \tilde{V} + M_v =: V_0.
\end{equation}}

The remaining proof on the boundedness of $\left\| D\left(\log_{x} x_{\pazocal{T}}\right)/dt \right\|_M$ follows an approach similar to that in \cite{li2025rendezvous}. First, using the Minkowski-coordinate representation of the logarithm map, we obtain
{\small \begin{equation}
    \log_x x_{\pazocal{T}} = \frac{d}{\sinh d}\,x_{\pazocal{T}} - \frac{d\cosh d}{\sinh d}\,x,
\end{equation}}
and differentiating with respect to $t$, we get
{\small \begin{equation}
    \frac{d}{dt}\!\left(\log_x x_{\pazocal{T}}\right) = I_1\, x_{\pazocal{T}} + I_2\, v_{\pazocal{T}} - I_3\, x - I_4\, v,
    \label{eq:log_derivative}
\end{equation}} where $I_i$'s for $i=1,2,3,4$ are defined as
{\small \begin{align*}
    I_1 &= \left(\frac{\sinh d - d\cosh d}{\sinh^2 d}\,\right)\dot{d},
    \quad
    I_2 = \left( \frac{d}{\sinh d}\right), \notag \\
    I_3 &= \left( \frac{\cosh d\sinh d - d}{\sinh^2 d}\ \right) \dot{d},
    \quad
    I_4 = {d\coth d}.
\end{align*}}
These coefficients are uniformly bounded as follows:
\begin{itemize}
\item  $I_2 \leq 1$, since $\sinh d \geq d$ for all $d \geq 0$.

\item $I_4 =d\coth d $ is continuous on $ [0, D_0]$ with $\lim_{d \to 0} I_4 =1$. Thus, $|I_4| \leq C_4 := \max_{s \in [0, D_0]} s\coth s < \infty$.

\item As $d \rightarrow 0$, the expression $\frac{\sinh d - d\cosh d}{\sinh^2 d} \to 0$. Thus, this fraction is continuous and bounded on $[0, D_0]$. Combining this result and bounds obtained in \eqref{eq:ddot_bound}, we get
    {\small \begin{equation*}
        |I_1| \leq C_1\,\tilde{V}, \quad
        C_1 := \max_{s \in [0,D_0]}\left|\frac{\sinh s - s\cosh s}{\sinh^2 s}\right|.
    \end{equation*}}

\item Similarly, $\frac{\cosh d \sinh d - d}{\sinh^2 d} \to 0$ as $d \to 0$, and is continuous and bounded on $[0, D_0]$. Thus,
   {\small  \begin{equation*}
        |I_3| \leq C_3\,\tilde{V}, \quad
        C_3 := \max_{s \in [0,D_0]}\left|\frac{\cosh s\sinh s - s}{\sinh^2 s}\right|.
    \end{equation*}}
\end{itemize}
The following bounds follow from the same argument as in \cite{li2025rendezvous}, applied to the agent-target pair:
{\small \begin{equation}\label{eq:cross1}
\begin{split}
|\langle v_{\pazocal{T}},\, x\rangle_M| &\leq \sinh(D_0)\,M_v, \\
|\langle v,\, x_{\pazocal{T}}\rangle_M| &\leq \sinh(D_0)\,V_0, \\
|\langle v_{\pazocal{T}},\, v\rangle_M| &\leq \cosh^2(D_0)\,M_v V_0.
\end{split}
\end{equation}}
Computing the covariant derivative $\log_x x_{\pazocal{T}}$ and taking its squared norm, we obtain
{\small 
\begin{align}\label{eq:cov_split}
\frac{D}{dt}\!\left(\log_x x_{\pazocal{T}}\right) = & \frac{d}{dt}\!\left[\log_x x_{\pazocal{T}}\right] + \left\langle \frac{d}{dt}\!\left[\log_x x_{\pazocal{T}}\right],\, x\right\rangle_M x \nonumber\\
\Rightarrow \left\|\frac{D}{dt}\log_x x_{\pazocal{T}}\right\|_M^2  = &\left\langle \frac{d}{dt}\log_x x_{\pazocal{T}},\, \frac{d}{dt}\log_x x_{\pazocal{T}}\right\rangle_M \nonumber \\
 &   + \left\langle \frac{d}{dt}\log_x x_{\pazocal{T}},\, x\right\rangle_M^2
\end{align}}
Expand the first term of \eqref{eq:cov_split} and using \eqref{eq:log_derivative}, we obtain
{\small \begin{equation} \label{eq:norm_bound}
\begin{split}
& \left\langle \frac{d}{dt}\log_x x_{\pazocal{T}},\, \frac{d}{dt}\log_x x_{\pazocal{T}}\right\rangle_M = -I_1^2 - 2I_1 I_3{(-\cosh d) }- 2I_1 I_4\langle x_{\pazocal{T}},\, v\rangle_M \\
& \qquad + I_2^2\|v_{\pazocal{T}}\|_M^2 - 2I_2 I_3\langle v_{\pazocal{T}},\, x\rangle_M - 2I_2 I_4\langle v_{\pazocal{T}},\, v\rangle_M - I_3^2 + I_4^2\|v\|_M^2
    \end{split}
\end{equation}}
 From \eqref{eq:bound:v}, \eqref{eq:cross1}, and using the bounds on $I_1, I_2, I_3$, and $ I_4$,
{\small \begin{equation}
\begin{split}
 & \left|\left\langle \frac{d}{dt}\log_x x_{\pazocal{T}},\, \frac{d}{dt}\log_x x_{\pazocal{T}}\right\rangle_M\right| \leq (C_1\tilde{V})^2 + 2C_1 C_3\tilde{V}^2\cosh(D_0)\\
& \qquad + 2C_1 C_4\tilde{V}\sinh(D_0)V_0 + M_v^2 + 2 C_3\tilde{V}\sinh(D_0)M_v \\
& \qquad + 2C_4\cosh^2(D_0)M_vV_0+ (C_3\tilde{V})^2 + C_4^2 V_0^2 =: B_1^2 < \infty
 \end{split}
\end{equation}}
Using \eqref{eq:log_derivative}, the second term in \eqref{eq:cov_split} becomes
{\small \begin{equation}\label{eq:bound:2nd}
    \left\langle \frac{d}{dt}\log_x x_{\pazocal{T}},\, x\right\rangle_M
    = -I_1\cosh d + I_2\langle v_{\pazocal{T}},\, x\rangle_M + I_3.
\end{equation}}
Using \eqref{eq:cross1} and upper bounds of $I_1, I_2, I_3$, we get from \eqref{eq:bound:2nd}
{\small \begin{equation}
    \left|\left\langle \frac{d}{dt}\log_x x_{\pazocal{T}},\, x\right\rangle_M\right|
    \leq C_1\tilde{V}\cosh(D_0) + \sinh(D_0)M_v + C_3\tilde{V} =: B_2 < \infty.
    \label{eq:proj_bound}
\end{equation}}
Substituting \eqref{eq:proj_bound} and \eqref{eq:norm_bound} into \eqref{eq:cov_split}, we obtain
{\small \begin{equation}
    \sup_{t \geq 0}\left\|\frac{D}{dt}\log_x x_{\pazocal{T}}(t)\right\|_M
    \leq \sqrt{B_1^2 + B_2^2} =: L < \infty.
    \label{eq:log_cov_bound}
\end{equation}}
 Thus, $\|D^2\tilde{v}/dt^2\|^2_M$ is bounded and $\lim_{t \rightarrow \infty}\tilde{v}=0$ (from the earlier part of the proof). Therefore, using the extension of Barbalat's lemma \ref{lemma:barbalat:hyperboloid},
$\lim _{t \rightarrow \infty} \frac{D\tilde{v}}{dt}(t) =0$.

From \eqref{eq:derv:vtilde}, this limit yields
{\small \begin{equation*}
\begin{split}
& \log_x x_{\pazocal{T}} = \frac{1}{k_p} \frac{D\tilde{v}}{dt}+\frac{k_d}{k_p}\tilde{v} \\ \Rightarrow & d(x(t),x_{\pazocal{T}}(t)) \leq \frac{1}{k_p} \left\| \frac{D\tilde{v}}{dt}\right\|_M+\frac{k_d}{k_p}\|\tilde{v}(t)\|_M
\end{split}
\end{equation*}}
Since $d(x(t),x_{\pazocal{T}}(t)) \ge 0$, we have 
{\small \begin{equation*}
\begin{split}
0 \leq \lim_{t\rightarrow \infty} & d(x(t),x_{\pazocal{T}}(t))\leq  \frac{1}{k_p} \lim_{t\rightarrow \infty} \left\| \frac{D\tilde{v}}{dt}\right\|_M + \frac{k_d}{k_p}\lim_{t\rightarrow \infty} \|\tilde{v}(t)\|_M \\
\Rightarrow 0 \leq \lim_{t\rightarrow \infty} & d(x(t),x_{\pazocal{T}}(t)) \leq \frac{1}{k_p}\cdot 0+\frac{k_d}{k_p}\cdot 0 = 0 \\
\implies \lim_{t\rightarrow \infty}&d(x(t),x_{\pazocal{T}}(t))=0
\end{split}
\end{equation*}}
Thus, this proof ensures both velocity and distance convergence with \eqref{eq:controller}. Hence, under the proposed controller \eqref{eq:controller}, an agent can track a moving target in $\mathbb{H}^n$.
\end{proof}


\section{ Leader-follower Formation Controller}\label{sec:leader-follower_formation}
We extend \eqref{eq:controller} to a formation control law, where agents with dynamics \eqref{eq:double_int_agent} form a desired geometry around the leader $L$ in a particular orientation. Let the leader's position is $x_L(t) \in \mathbb{H}^n$. We define a virtual target $i^*$ for the $i$th agent to track, where its position $x_i^*(t) \in \mathbb{H}^n$ around the leader in a particular orientation is computed using exponential map (see Definition \ref{def:log:exp}) as
{\small \begin{equation}
    \label{eq:virtual:target}
    x_i^*(t)=\operatorname{exp}_{x_L(t)}(\xi_i^*(t)).
\end{equation}}
{Note that $\xi_i^*(t)\in T_{x_L(t)}\mathbb{H}^n$ is the tangent vector to determine the desired relative position and orientation of $i^*$ with respect to $L$.} To preserve the prescribed formation around the moving leader, $\xi_i^*$ is parallel transported along the leader's trajectory, that is {$\frac{D}{dt}\xi_i^* \equiv 0$}. Thus, $i^*$ always maintains a fixed relative configuration with respect to $L$. 

By definition of parallel transport \cite{li2025rendezvous}, $\xi^*_i(t)$ in terms of its initial configuration $\xi^*_i(0)$ is expressed as
{\small
\begin{equation}\label{eq:position:virtual:target}
    \xi^*_i(t)=P_{x_L(0)}^{x_L(t)} \xi^*_i(0).
\end{equation}}
where $P$ denotes the parallel transport operator defined in Section \ref{sec:prelim}). Since $\frac{D}{dt}\xi_i^*\equiv0$, using the metric compatibility of the Levi-Civita connection, we get
\begin{center}
$ \frac{d}{dt} \|\xi_i^*\|_M^2 = \frac{d}{dt} \langle \xi_i^*, \xi_i^* \rangle_M = 2\langle \frac{D}{dt} \xi_i^*, \xi^*_i \rangle_M=0 $.
\end{center}
Thus, denoting $\|\xi_i^*\|_M^2 = r_i$, we write $\xi_i^*(t)= r_i q(t)$, where $q(t)$ is the unit vector in $T_{x_L(t)}\mathbb{H}^n$ along $\xi^*_i(t)$. From \eqref{eq:virtual:target},
{\small
    \begin{equation}\label{eq:pos:virtual:target}
        x_i^*(t)= \operatorname{exp}_{x_L(t)}(\xi_i^*(t))= \cosh(r_i)x_L(t)+\sinh(r_i){q_i(t)}.
\end{equation}}    
Writing the intrinsic covariant derivative as the tangential projection of the ambient derivative of $q(t)$, we get
{\small
\begin{equation} \label{eq:virtual:target:derv}
     \frac{D}{dt}q(t)= \dot{q}(t)+ \langle \dot{q}(t), x_L(t)\rangle_Mx_L(t).
\end{equation}
}
Since $q(t)$ is a unit tangent vector of $x_L(t)$, differentiating $\langle x_L(t), q(t) \rangle_M \equiv 0$, we get $ \langle \dot{q}(t), x_L(t) \rangle_M = - \langle v_L(t), q(t) \rangle_M $. Thus, \eqref{eq:virtual:target:derv} becomes
{ \small \begin{equation}\label{eq:virtual:target:derv:next}
    \frac{D}{dt}q(t)=  \dot{q}(t)- \langle v_L(t), q(t)\rangle_Mx_L(t).
\end{equation}}
Furthermore, since $ \frac{D}{dt} \xi_i^* \equiv0$ and $r_i$ is constant, $\frac{D}{dt} q(t) \equiv 0$. Thus, 
\eqref{eq:virtual:target:derv:next} gives $\dot{q}(t)= \langle v_L(t), q(t) \rangle_M x_L(t)$. Therefore, velocity of $i^*$ is calculated as
{\small \begin{equation}\label{eq:vel:virtual:target}
\begin{split}
    v_i^*(t)= \dot{x}_i^*(t) & = \cosh (r_i) v_L(t) + \sinh (r_i) \dot{q}(t)\\\
    & = \cosh (r_i) v_L(t) + \sinh (r_i) \langle v_L(t), q(t) \rangle_M x_L(t)
\end{split} 
\end{equation}}
{The following corollary establishes extension of \eqref{eq:controller} to a leader-follower protocol.

\begin{corollary}
Consider $N$ agents $i=1,2 ,\ldots, N$ with motion \eqref{eq:double_int_agent} and a moving leader satisfying Assumption \ref{ass:target:speed:bound} such that $\|v_L\|_M \leq M_L < \infty$. Let the tracking controller \eqref{eq:controller} is modified with target position $x_i^*(t)$ in \eqref{eq:pos:virtual:target} and velocity $v_i^*(t)$ in \eqref{eq:vel:virtual:target}. Then the following hold for $i=1, 2, \ldots, N$:
\begin{enumerate}[(i)]
    \item The velocity of $i^*$ is bounded as
    {\small
\begin{equation*}
\|v_i^*(t)\|_M^2 \leq \cosh(r_i)^2\|v_L(t)\|^2_M  \leq \cosh(r_i)^2{M_L^2}
\end{equation*}
}
\item Each agent $i$ tracks the trajectory of $i^*$ on $\mathbb{H}^n$ around a moving leader $L$ ensuring
{\small
\begin{equation*}
 \lim _{t \rightarrow \infty} \left\| \log_{x_i(t)} x_{i}^*(t) \right\|_M^2 = 0 \qquad \lim _{t \rightarrow \infty} \left\| v_i(t) - P_{x_{i}^*(t)}^{x_i(t)} v_{i}^*(t)\right\|_M^2 =0
\end{equation*}}
This ensures a leader-follower formation around a moving leader $L$ on $\mathbb{H}^n$. 
\end{enumerate}
\end{corollary}

\begin{proof}
\textbf{Part 1:} From \eqref{eq:vel:virtual:target}, we get the bound of $v^*_i$ as 
{\small
\begin{equation*}
\begin{split}
      \|v_i^*(t)\|_M^2 &= \cosh(r_i)^2\|v_L(t)\|^2_M -\sinh(r_i)^2 \langle v_L(t), q(t) \rangle_M^2\\
      & \leq \cosh(r_i)^2\|v_L(t)\|^2_M  \leq \cosh(r_i)^2{M_L^2}.
\end{split}
\end{equation*}
}
Thus, $i^*$ satisfies Assumption \ref{ass:target:speed:bound}.

\textbf{Part 2:} After substituting $x_i^*$ and $v_i^*$ in stead $x_{\pazocal{T}}$ and $v_{\pazocal{T}}$ in \eqref{eq:controller}, the modified controller ensure tracking of $i^*$ by $i$. The proof is same as that of Theorem \ref{thm:tracking}. Hence, a leader-follower formation is formed around a moving leader $L$ on $\mathbb{H}^n$. 
\end{proof}
}


\section{Numerical Simulations}\label{sec:sim}
In this section, we demonstrate the effectiveness of the control protocol \eqref{eq:controller} in {two scenarios} on the hyperboloid $\mathbb{H}^2$.
\subsection{Target-tracking with Single Agent}\label{subsec:results:target:track}
In this example, two different target trajectories cases are considered: geodesic and spiral. For geodesic trajectory, the target $\pazocal{T}$ moves at a constant speed of $v_r=0.25$, whereas $\pazocal{T}$ to execute the spiral trajectory evolves along a hyperbolic spiral on the hyperboloid, with position given by $\mathbf{x_{\pazocal{T}}}(t)=\cosh(r_0+v_rt)\mathbf{p}+\sinh(r_0+v_rt)\mathbf{u}_\theta(t)$, where $\mathbf{u}_\theta(t)=\cos(\omega t)\,\mathbf{e}_1+\sin(\omega t)\,\mathbf{e}_2$, $\mathbf{p}=(1,0,0)^\top$ is the base point, $\mathbf{e}_1=(0,1,0)^\top$ and $\mathbf{e}_2=(0,0,1)^\top$ are the tangent basis vectors at $\mathbf{p}$, and $r_0=0.3, v_r=0.2,{\omega=\frac{\omega_0}{\cosh{r_0+v_rt}}, \omega_0=0.2}$. {Note that for the geodesic trajectory, $\|v_\pazocal{T}\|_M$ is constant, and for the spiral trajectory, $\footnotesize{\|v_\pazocal{T}(t)\|_M = \sqrt{v_r^2+ \omega_0^2\tanh ^2(r_0+v_rt)}< \sqrt{v_r^2+ \omega_0^2}}$. Hence, $\|v_{\pazocal{T}}(t)\|_M$ is bounded for both geodesic and spiral trajectories.} Initial positions of target and agent in the geodesic and spiral scenarios are set as given in Table \ref{tab:initial:target:track}. The initial speed of the agent is set to zero. {Parameters $k_p$ and $k_d$ are set to $8$ and $4$, respectively.} It can be verified that initial positions of agent and target lie on the hyperboloid $\mathbb{H}^2$.

\begin{table}[h!]
    \centering
    \scalebox{0.8}{
    \begin{tabular}{|c|c|c|}
    \hline
         &Spiral&Geodesic \\
         \hline
         Target& $(1.045,0.304,0)$& $(1,0,0)$\\
         \hline
         Agent & $(1.810,-0.628,1.372)$&$(1.810,-0.628,1.372)$\\
         \hline
    \end{tabular}}
    \caption{Initial Conditions for Target Tracking}
    \label{tab:initial:target:track}
    \vskip -10pt
\end{table}

\begin{figure}[h!]
    \centering
    \subfloat[For geodesic trajectory on $\mathbb{H}^2$]{\label{fig:tracking:geodesic:trajectory}
    {\includegraphics[width=4cm, height = 3.75cm]{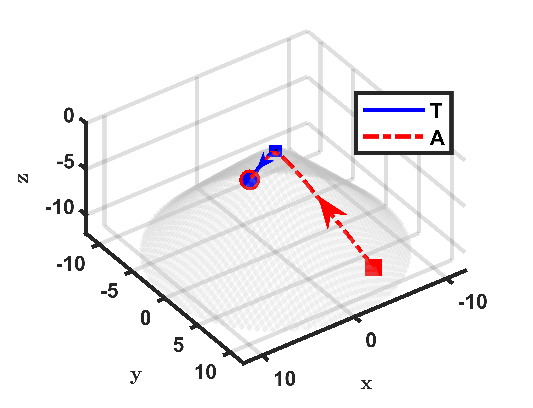}}}%
    \quad
    \subfloat[For spiral trajectory on $\mathbb{H}^2$]{\label{fig:tracking:spiral:trajectory}
    {\includegraphics[width=4cm, height = 3.75cm]{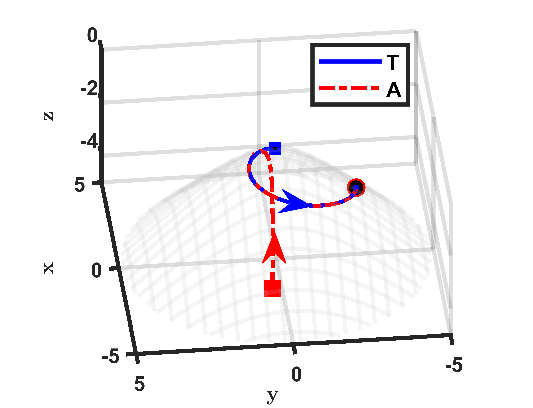}}}%
    \caption{Trajectories of agent and target $\pazocal{T}$ ensuring target-tracking under \eqref{eq:controller}. In two examples, $\pazocal{T}$ is executing geodesic (Fig. \ref{fig:tracking:geodesic:trajectory}) and spiral (Fig. \ref{fig:tracking:spiral:trajectory}) trajectories on $\mathbb{H}^2$}\label{fig:tracking:traj}
    \vskip -15pt
\end{figure}

\begin{figure}[h!]
    \centering

\subfloat[Distance errors over time]{\label{fig:tracking:error:distance}
{\includegraphics[width=4.1cm, height = 2.85cm]{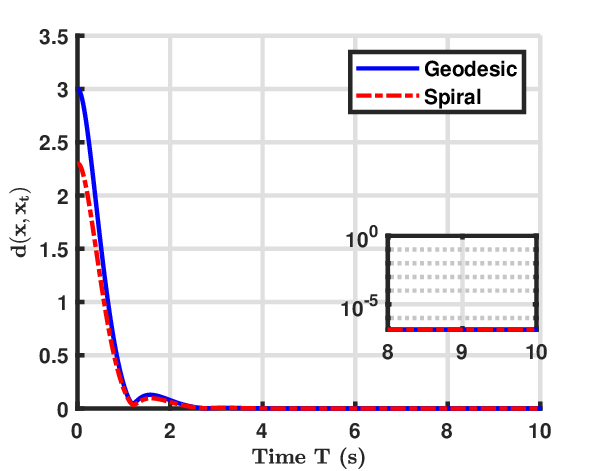}}}%
~
\subfloat[Speed errors over time]{\label{fig:tracking:error:speed}
{\includegraphics[width=4.1cm, height = 2.85cm]{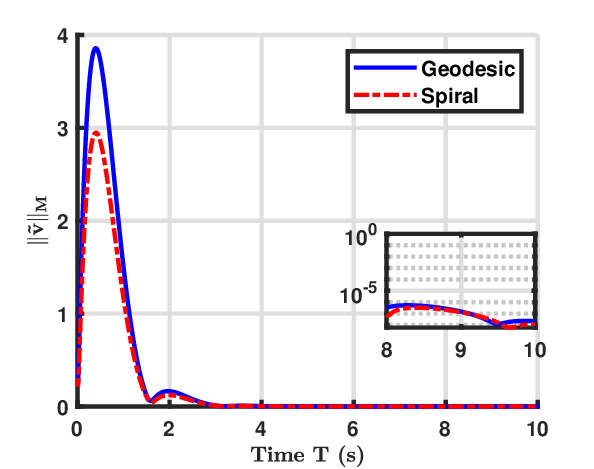}}}%
\caption{Evolution of position (Fig. \ref{fig:tracking:error:distance}) and velocity (Fig. \ref{fig:tracking:error:speed}) errors while $\pazocal{T}$ is executing geodesic and spiral trajectories on $\mathbb{H}^2$}\label{fig:tracking:error}
\vskip -10pt
\end{figure}

The performance of \eqref{eq:controller} is demonstrated for $\pazocal{T}$ moving along a geodesic (Fig. \ref{fig:tracking:geodesic:trajectory}) and a spiral (Fig.\ref{fig:tracking:spiral:trajectory}) trajectory in $\mathbb{H}^2$. The agent is initialized at a position (red square) distinct from the target's (blue square), and converges to the target under \eqref{eq:controller} reaching position and velocity agreement at endpoints (marked by circles). This is further confirmed by asymptotic convergence of distance and speed errors as shown in Fig. \ref{fig:tracking:error}. {Thus, the simulations confirm that agent successfully tracks a moving target on $\mathbb{H}^2$ under \eqref{eq:controller}.}

\subsection{Multi-agent Formation on Hyperboloid $\mathbb{H}^2$}\label{subsec:results:formation:control}
For the formation-control example, we consider a {three agent}  MAS with {a leader $L$} following the same trajectories (geodesic and spiral) with same initial condition as the target $\pazocal{T}$ in Section \ref{subsec:results:target:track}.{ Agents' initial positions and orientations are listed in Table \ref{tab:initial:formation:control}. The virtual targets $i^*$ are placed at equiangular spacing and at equal distances of $0.8$ units from $L$}. {Parameters $k_p$ and $k_d$ are set to $4$.} Note that the initial positions lie on $\mathbb{H}^2$. Similar to Section \ref{subsec:results:target:track}, each agent starts from the rest.
\begin{table}[t]
    \centering
    \scalebox{0.85}{
    \begin{tabular}{|l|c|c|}
        \hline
        \multirow{2}{*}{} & \multicolumn{2}{c|}{\bf Spiral} \\
        \cline{2-3} 
        Agent & Position & Orientation \\
        \hline
        1 & $(1.9535,1.6361,-0.3731)$ & $(0.2436,0.8362,0)$ \\
        \hline
        2 & $(1.0836,0.3559,-0.2178)$ & $(-0.1218,-0.4181,0.6928)$ \\
        \hline
        3 & $(1.3945,-0.3644,-0.9009)$ & $(-0.1218,-0.4181,-0.6928)$ \\
        \hline
        \multirow{2}{*}{} & \multicolumn{2}{c|}{\bf Geodesic} \\
        \cline{2-3} 
        & Position & Orientation \\
        \hline
        1 &$(1.9648,1.6870,0.1202)$ & $(0,0.8,0)$ \\
        \hline
        2 & $(1.4967,-1.10133,0.1652)$ & $(0,-0.4,0.6928)$ \\
        \hline
        3 & $(1.5161,0.6028,-0.9671)$ & $(0,-0.4,-0.6928)$ \\
        \hline
    \end{tabular}
    }
    \caption{Initial conditions for formation control}
    \label{tab:initial:formation:control}
    \vskip -15pt
\end{table}

\begin{figure}[h!]
    \centering
    \subfloat[For Geodesic trajectory of $L$ in $\mathbb{H}^2$]{\label{fig:formation:geodesic:trajectory}
    {\includegraphics[width=4.2cm, height = 3.75cm]{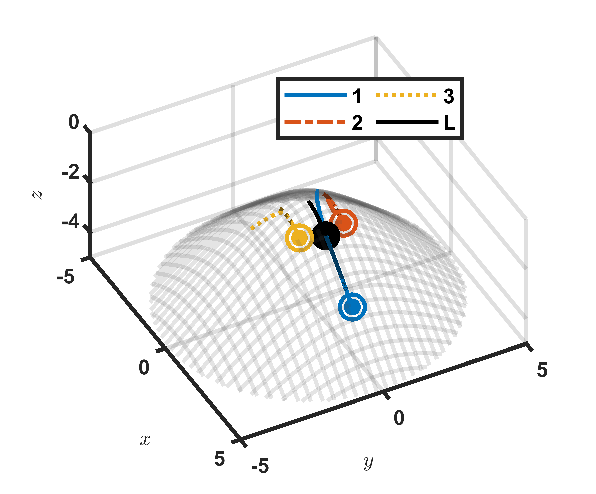}}}%
    ~
    \subfloat[For Spiral trajectory of $L$ in $\mathbb{H}^2$]{\label{fig:formation:spiral:trajectory}
    {\includegraphics[width=4cm, height = 3.75cm]{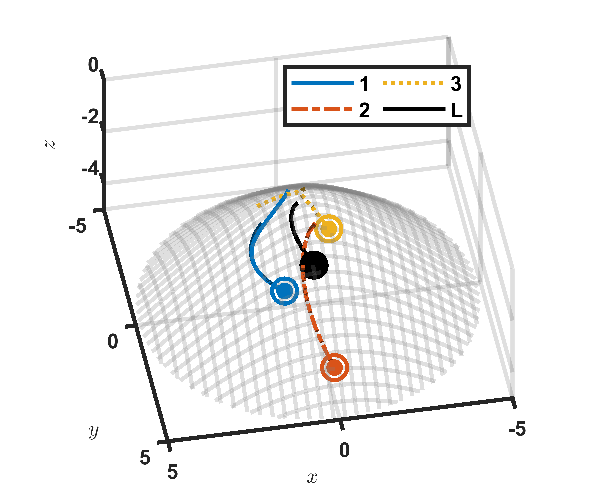}}}%
    \quad
    \subfloat[For Geodesic trajectory of $L$ in $\mathbb{H}^2$ ($X-Y$ projection)]{\label{fig:formation:geodesic:2D}
    {\includegraphics[width=4cm, height = 2.95cm]{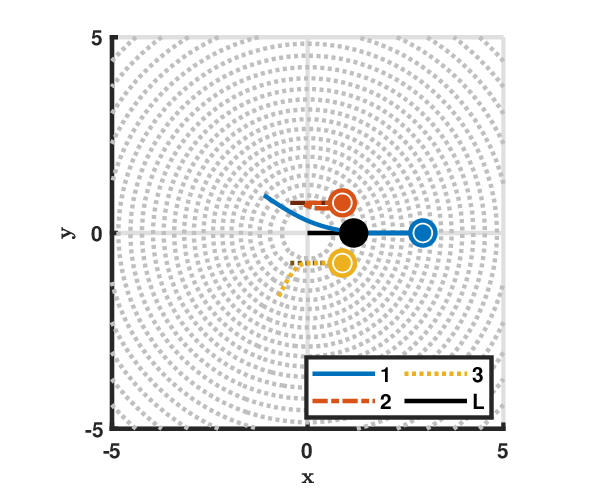}}}%
    ~
    \subfloat[For spiral trajectory of $L$ in $\mathbb{H}^2$ ($X-Y$ projection)]{\label{fig:formation:spiral:2D}
    {\includegraphics[width=4cm, height = 2.95cm]{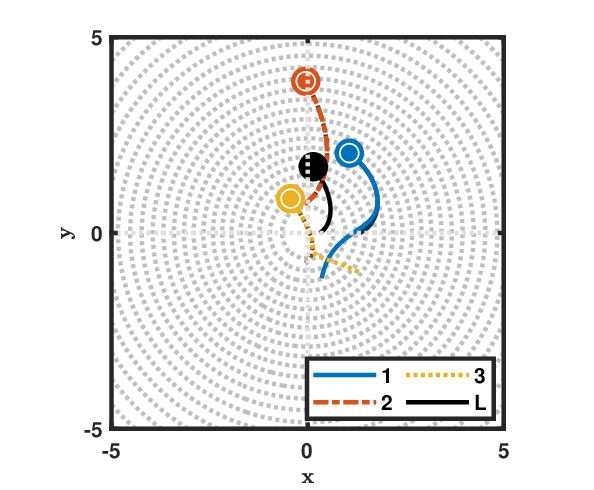}}}%
    \caption{Trajectories of followers and leader $L$ showing formation in $\mathbb{H}^2$. In two examples, $L$ is executing geodesic (Fig. \ref{fig:formation:geodesic:trajectory} and Fig. \ref{fig:formation:geodesic:2D} ($X-Y$ projection)) and spiral (Fig. \ref{fig:formation:spiral:trajectory} and Fig. \ref{fig:formation:spiral:2D} ($X-Y$ projection)) trajectories in $\mathbb{H}^2$}\label{fig:formation}
    \vskip -10pt
\end{figure}

\begin{figure}[h!]
    \centering
    \subfloat[For geodesic trajectory of $L$]{\label{fig:formation:agentdistance:geodesic}
    {\includegraphics[width=4cm, height = 2.85cm]{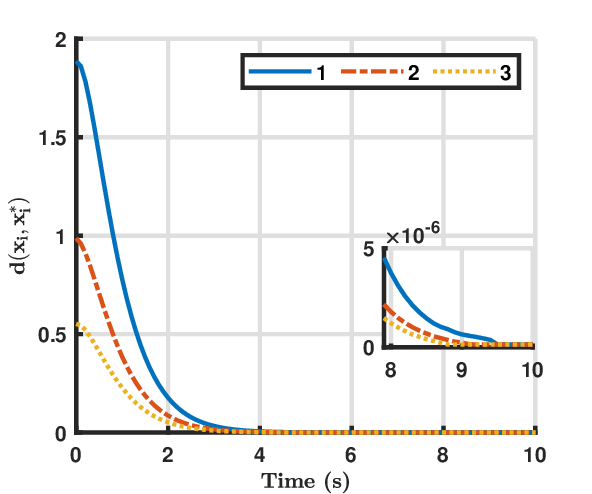}}}~
    \subfloat[For spiral trajectory of $L$]{\label{fig:formation:agentdistance:spiral}
    {\includegraphics[width=4cm, height = 2.85cm]{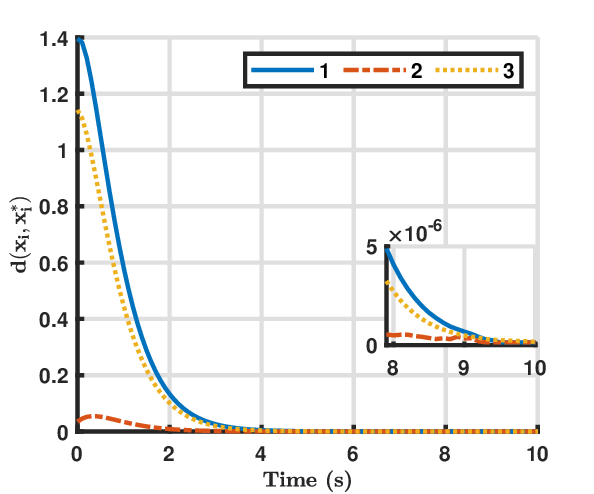}}}%
    \caption{Convergence of leader-follower distances in $\mathbb{H}^2$ ensuring desired formation around a leader $L$ moving with bounded velocity}\label{fig:formation:agentdistance}
    \vskip -15pt
\end{figure}

{The trajectories of three agents converging to the desired formation in $\mathbb{H}^2$ are shown in Fig. \ref{fig:formation}. Initialized at different positions, each agent converges to an equidistant spacing around $L$. The achievement of desired formation is further confirmed in Fig. \ref{fig:formation:agentdistance} demonstrating that each agent asymptotically reaches at a distance of $0.8$ units from $L$ on $\mathbb{H}^2$.}



\section{Conclusion}\label{sec:conc}

This paper proposes a novel controller for a double-integrator agent to track a moving target on $\mathbb{H}^n$. The control law is designed using the intrinsic Riemannian structure of the hyperboloid, which includes the logarithm map, covariant derivative, and parallel transport. Thus, all the terms in the proposed controller are coordinate-free, easing its implementation in the ambient Minkowski space without any additional chart computations. Using Barbalat's lemma for the hyperboloid, we have proved the asymptotic convergence of the tracking error for a target with a possibly unbounded trajectory and bounded velocity. The proposed controller is later extended to a leader-follower formation control, where the agents attain a desired geometry in $\mathbb{H}^n$ around a moving leader. The efficacy of the proposed controllers is illustrated with numerical simulations considering geodesic and spiral target trajectories. {The present framework also opens several directions for further investigations. In particular, the feedforward component of the controller relies on knowledge of the target's acceleration, while the convergence analysis assumes exact state information and ideal continuous-time feedback. These requirements motivate the development of acceleration-free control strategies and the investigation of robustness under disturbances, measurement noise, communication delays, and actuator constraints. Extending the present stability analysis to accommodate such practical considerations would help bridge the gap between the geometric guarantees established here and their implementation in realistic multi-agent systems on hyperbolic spaces.}

\section*{Tool use} No AI tool was used either for the conceptualization or for the exposition of this manuscript. 

\bibliography{main}
\bibliographystyle{IEEEtran}

\end{document}